\documentclass[11 pt]{amsart}
\usepackage{amscd,amsfonts,amssymb,amsmath}
\usepackage{hyperref}
\usepackage{epsfig}
\usepackage{mathtools}
\usepackage{tabu}
\usepackage{tikz-cd}
\newtheorem{theorem}{Theorem}[section]

\theoremstyle{definition}

\newtheorem{problem}[theorem]{Problem}
\newtheorem{definition}[theorem]{Definition}

\newtheorem{remark}[theorem]{Remark}
\numberwithin{equation}{subsection}
\newtheorem*{ack}{Acknowledgement}
\usepackage[all,cmtip]{xy}

\usepackage{graphicx} 

\begin{document}

\title{On certain properties of pseudo-Riemannian Bertrand manifolds}
\author{Denis A. Fedoseev}

\date{}

\address{Moscow state university}
\email{denfedex@yandex.ru}


\begin{abstract}
In the present article we prove and discuss several properties of pseudo-Riemannian Bertrand manifolds, and give an overview of the state of the theory together with problems for future work. In particular, we prove that there are no pseudo-Riemannian completely Bertrand systems --- dynamical systems of movement in a central field on a pseudo-Riemannian 2-dimensional manifold of revolution such that all trajectories of movement are closed.
\end{abstract}
\maketitle

\bigskip

\section{A brief history of Bertrand problem}\label{introduction}
The terms ``Bertrand manifold'' and ``Bertrand problem'' originate from the 1873 paper~\cite{Bertrand} by the prominent mathematician J.~Betrand which was devoted to the following mechanical problem: {\em consider a point moving according in a central potential field in a 3-dimensional Euclidean space; provided that all trajectories of the point (with initial velocity not greater than a certain limit) are closed, what can be said about the central potential in question?}

The problem was successfully solved by Bertrand himself (although, under certain additional assumption which were not explicitly in the original article and were recovered much later when the paper was carefully scrutinised, see~\cite{ZKF}). He proved that there are exactly two suitable potentials: Newton's gravity potential and oscillator potential.

In the following century and a half the problem was repeatedly reformulated and generalized~\cite{Liebmann, Darboux1, Darboux2, Perlick, Santoprete, ZKF}. Most prominent classical results in this area are due to Darboux~\cite{Darboux1, Darboux2}. In these works Darboux considered an arbitrary manifold of revolution (implicitly assuming that it does not have {\em equators}, see below) and asked roughly the same question: {\em if all trajectories of a certain class defined as a subset of possible initial conditions of the trajectories, are closed, what can be said on the central potential of the system?} His solution was later refined by Perlick~\cite{Perlick} and Santoprete~\cite{Santoprete}, and the final form of the generalized Bertrand problem was given (and solved for manifolds without equators) in~\cite{ZKF}. 

The last mentioned paper contained the detailed analysis of the previous approaches to the problem, classified its formulations and contained the most complete solution to date. To be precise, in~\cite{ZKF} five types of Bertrand potential were introduced: {\em closing, locally closing, semilocally closing, weakly closing} and {\em strongly closing}. It was proved that for the case of manifolds of revolution without equators those five classes coincide, and an explicit form of Riemannian metric was presented which allowed the existence of the potentials of those classes. Furthermore, the manifolds (called Bertrand manifolds of the corresponding type --- Bertran, locally Bertrand, semilocally Bertrand, etc.) were classified by three parameters, and the explicit number and form of the potentials were found depending on the parameters. It was shown that there are at most two potentials (in which case they are generalizations of Newton and oscillator potentials), and for a certain subfamily of manifolds there exists only one (oscillator) potential.

Note that these results were obtained only under assumption that the configurational manifold of the problem does not have any equators. At the same time, it was clear that there exist bertrand manifolds with equators. For example, the so-called {\em Tannery surfaces} --- surfaces of revolution, on which all geodesics are closed, see~\cite{Besse}, --- are clearly Bertrand manifolds which admit constant Bertrand (and even {\em completely Bertrand}) potential. Those manifolds are well-studied and do have equators. Hence, there exist Bertrand manifolds with equators.

Bertrand problem on Riemannian maniflds of revolution with equators was subsequently tackled in~\cite{KudrFed1, KudrFed2, KudrFed3} and completely solved for the cases of {\em completely, stably, weakly} and {\em strongly} Bertrand systems, as well as {\em Bertrand} systems. Together with the results on the case of manifolds without equators and the intersection diagram of the Bertrand classes (see~\cite{KudrFed3}) the only case which remained unsolved as of now is the case of locally Bertrand systems on manifolds of revolution with equators.

One should notice that the results listed above were obtained for the case of {\em Riemannian} manifolds of revolution. Nevertheless, an investigation of the pseudo-Riemannian case was carried out by Zagryadskii~\cite{Zagr}. The problem was solved for the case of manifolds without equators (parallel to the results of~\cite{ZKF}). Somewhat surprisingly, the answer turned out to strongly resemble the Riemannian case: Bertrand manifolds were described by a subfamily of the same 3-parameter family as in the Riemannian case (though, in constrast with the latter, not all values of parameters give a pseudo-Riemannian Bertrand manifolds, and the radial coordinate interval differs in the two cases).

The case of manifolds with equators remained completely unstudied. In the present article we provide some steps towards the solution of this problem. In particular, we prove that there are no {\em completely Bertrand systems} in the pseudo-Riemannian case (either with or without equators) and describe the general form of locally Bertrand pseudo-Riemannian manifolds in the case when strongly stable orbits do not coincide with the equator(s).

\section{Preliminaries}
\subsection{Necessary definitions}

In the present paper we consider mechanical systems on manifolds of revolution, that is a two-dimensional noncompact manifold
$S \simeq (a, b) \times S^1 $ endowed with a (pseudo) Riemannian metric, which in the natural ``polar'' coordinates $(r,\varphi)$ can be written in the form
$$ds^2 = dr^2 + \epsilon f^2(r)d\varphi^2,$$
where $f(r)$ is a smooth (or analytic) function on $(a,b)$ that has no zeros on this interval, and $\epsilon\in\{+1, -1\}$ is a sign which makes the metric either Riemannian or pseudo-Riemannian. The {\em equator} of a surface of revolution is such $r_0\in (a, b)$ that $f'(r_0)=$. Alternatively, we should call the whole parallel $\{r_0\}\times S^1$ the equator of the manifold.

Given a manifold of revolution $S$, we consider a natural mechanical system with two degrees of freedom with the configuration manifold $(S, ds^2)$ and the central potential $V = V (r)$. It is known that the {\em energy} $H = H(x,\dot{x }) = \frac{1}{2}(\dot{r}^2 + \epsilon f^2(r)\dot{\varphi}^2) + V(r)$ and the {\em kinetic moment} $K = K(x,\dot{x}) = \epsilon f^2(r)\dot{\varphi}$ are constant along any phase trajectory $(x(t),\dot{x}(t))$ of the given system. The function $W(r) = V(r)+\epsilon \frac{K^2}{2f^2(r)}$ is called the {\em effective potential} of the system.

\begin{definition}[\cite{KudrFed3}]
Let us call a {\em trajectory} a solution $x(t) = (r(t),\varphi(t))$ of the equations of motion defined on the maximal interval $(t_0,t_1) \subset \mathbb{R}^1$, and the {\em orbit} is the image $\{x(t) | t \in (t_0,t_1)\} \subset S$ of the mapping $x(t)$.
\begin{enumerate}
\item The orbit of a particle moving along the surface $S$ according to the law of force given by a central potential $V (r)$ (i.e., a potential depending only on $r$) is said to be {\em circular} if it coincides with the orbit of action of the rotation group. A trajectory is called circular if the corresponding orbit is circular.
\item A trajectory is said to be {\em bounded} if it is defined on the whole time axis $t \in \mathbb{R}^1$ and its image is contained in a certain compact set $[r_1, r_2] \times S^1 \subset (a, b) \times S^1$. An orbit is called bounded if the corresponding trajectory is bounded.
\item A trajectory (and the corresponding orbit) is said to be {\em singular} (see~\cite{Santoprete}) if the value of the kinetic moment integral K on this trajectory is zero, i.e., $\varphi = const$.
\end{enumerate}
\end{definition}

Now let us define several types of central potentials and the corresponding mechanical systems following~\cite{KudrFed3}.

\begin{definition}
\begin{enumerate}
	\item A central potential $V (r)$ on a surface $S$ is called a {\em closing potential} if
	\begin{itemize}
		\item[($\exists$)] there exists a nonsingular, bounded, noncircular orbit $\gamma$ in $S$;
		\item[($\forall$)] each nonsingular bounded orbit in S is closed.
	\end{itemize}
	\item A potential $V (r)$ is said to be {\em locally closing} if
	\begin{itemize}
		\item[$(\exists)^{loc}$] there exists a strongly stable circular orbit $\{r_0\} \times S^1$ in $S$;
		\item[$(\forall)^{loc}$] for any strongly stable circular orbit $\{r_0\} \times S^1$ in $S$, there exists a number $\varepsilon > 0$ such that each nonsingular bounded orbit entirely contained in the ring $[r_0 - \varepsilon, r_0 + \varepsilon] \times S^1$ and having a kinetic moment level in the interval $(K_0 - \varepsilon, K_0 + \varepsilon)$ is closed; here $K_0$ is the value of the kinetic moment on the corresponding circular trajectory.
	\end{itemize}
	\item A potential $V(r)$ is called {\em completely closing} if all nonsingular orbits are closed. If all nonsingular geodesics are closed (i.e., the constant potential $V_0 = const$ is completely closed), then the manifold $(S,ds^2)$ is called a {\em Tannery manifold}.
\end{enumerate}
\end{definition}

\begin{definition}
	Assume that a potential $V(r)$ belongs to one of the classes defined above. In this case, we say that the manifold $(S, ds^2)$, the triple $(S, ds^22, V)$, and the corresponding natural mechanical system have Bertrand type. Namely, they are said to be Bertrand, locally Bertrand, or completely Bertrand if the potential $V(r)$ is a closing, locally closing, or completely closing potential, respectively.
\end{definition}

\begin{remark}
	Note than unlike the Riemannian case, geodesics on pseudo-Riemannian manifolds form three distinct classes: spacelike, timelike and lightlike. Therefore, speaking of Tannery surfaces (and, by extension, of Bertrand manifolds) it is reasonable to differentiate spacelike, timelike and lightlike Tannery (Bertrand) manifolds. In the present article, if no adjective is used, we shall presume that the property of being closed is satisfied by {\em all} geodesic (trajectories).
\end{remark}

\subsection{Known classification theorems}

In the present section we recall some of the know classification theorems for Riemannian and pseudo-Riemannian Bertrand manifolds with and without equators. We also mention several theorems on Tannery surfaces.

In case of manifolds without equators, both Riemannian and pseudo-Riemannian case may be summirised in the following theorem (see~\cite{ZKF} for Riemannian and~\cite{Zagr} for pseudo-Riemannian parts of the theorem).

\begin{theorem} \label{thm:bert}
	Let $S$ be a (pseudo) Riemannian manifold of revolution without equators. We consider it to have the metric $ds^2=dr^2+\epsilon f^2(r)d\varphi^2$ in polar coordinates $(r, \varphi \mod 2\pi)$. Then the following statements hold:
	\begin{enumerate}
		\item classes of closing, locally, semilocally, weakly, and strongly closing potentials coincide;
		\item there exists at least one potential of those classes if and only if there exists such coordinate change $\theta = \theta(r)$ that in coordinates $(\theta, \varphi \mod 2\pi)$ the metric takes form $$ds^2=\frac{d\theta^2}{(\theta^2+c-\delta\theta^{-2})^2}+\frac{d\varphi^2}{\mu(\theta^2+c-\delta\theta^{-2})},$$ where $\mu\in \mathbb{Q}, c,\delta\in \mathbb{R}$ are constants, and the interval for $\theta$ is non-empty (see below);
		\item the intervals for $\theta$ depend on the case (Rimannian or pseudo-Riemannian) and on the parameters $c,\delta$; they are summarized in the table below:
		\begin{center}
		\begin{tabular}{|c|c|c|}
		\hline
		 Value of $(c,\delta)$ & Possible values for $\theta$& Possible values for $\theta$ \\ 
		   & in Riemannian case & in pseudo-Riemannian case \\ 
		\hline
		 $c\ge 0, \delta=0$ & $(0,\infty)$ & $\emptyset$ \\  
		\hline
		 $c<0, \delta=0$ & $(\sqrt{-c}, \infty)$ & $(0, \sqrt{-c})$ \\ 
		\hline
		 $\delta>0$ & $(\theta_2,\infty)$ & $(0,\theta_2)$ \\
		\hline
		 $\delta<0, c<0, c^2+4\delta>0$ & $(0,\theta_1) \cup (\theta_2,\infty)$ & $(\theta_1,\sqrt[4]{-\delta})\cup(\sqrt[4]{-\delta},\theta_2)$ \\
		\hline
		 $\delta<0, c>0, c^2 +4\delta>0$ or $c^2 + 4\delta \le 0$ & $(0,\sqrt[4]{-\delta})\cup(\sqrt[4]{-\delta},\infty)$ & $\emptyset$ \\
		\hline
		\end{tabular}
		\end{center}
		where $\theta_1, \theta_2$ are roots of the denominator.
		\item there exist exactly two potentials of the above-mentioned classes if and only if $\delta=0$ (and the corresponding interval for $\theta$ is non-empty), those potentials are generalized gravitational and generalized oscillator ones; for other values of $c$ and $\delta$ (provided that the corresponding interval for $\theta$ is non-empty) there exists exactly one potential (generalized oscillator).
	\end{enumerate}
\end{theorem}

\begin{remark}
	In the last two cases the value $\sqrt[4]{-t}$ corresponds to an equator which is removed from the manifold.
\end{remark}

Here we see that pseudo-Riemannian manifolds in a way ``compliment'' the Riemannian ones: the intervals for $\theta$ in pseudo-Riemannian case are complimentary to Riemannian ones in the set $\mathbb{R}\setminus\{\text{equators}\}$. This result is pretty natural if we take into account that the general form of the metric in both cases is the same, hence for each case $\theta$ should take such values that the sign of the second summand is appropriate.

Now we formulate the classification theorem for completely Bertrand manifolds (without the condition on equators), see~\cite{KudrFed3}.

\begin{theorem}[\cite{KudrFed3}] \label{thm:comp_bert_rim}
	Completely Bertrand Riemannian manifolds of revolution $(S,ds^2)$ together with completely closing central potentials $V$ on $S$ form, up to conjugacy of triples $(S,ds^2,V)$, five families (i)--(v) from the following list:
	\begin{enumerate}
	\item[(i)] ``rational'' (i.e., whose complete angle at the pole is a rational multiple of $2\pi$), spherical (i.e., having two poles and one equator), ``spheroidal'' (i.e., of a constant positive curvature) two-dimensional Riemannian manifolds of revolution (that form a two-parameter family containing a one-parameter family of round spheres with punctured poles) with the gravitational potential on them;
	\item[(ii)] ``rational'' cones with flat metrics (forming a one-parameter family containing the punctured Euclidean plane) with the oscillator potential on them;
	\item[(iii)] ``northern hemispheres'' of manifolds from item (i) with the oscillator potential on them;
	\item[(iv)] ``rational'' pear-shaped two-dimensional Riemannian manifolds of revolution from~\cite{ZKF} (forming a three-parameter family) with the oscillator potential on them, where the sign of the potential is such that the potential increases from the ``main'' pole to the ``additional,'' i.e., the ``main'' pole is attractive whereas the ``additional'' pole is repulsive;
	\item[(v)] all Tannery manifolds classified in~\cite{Besse} that form a ``functionally-two-parameter'' family (whose parameters are two positive real numbers, one of which is rational, and an odd smooth function $h\colon (-1, 1) \to (-1, 1)$ (defined up to the substitution $h \mapsto -h$) and that are spherical and spheroidal $(h \equiv 0)$ or pear-shaped $(h\neq 0)$ with a constant potential on them.
	\end{enumerate}
\end{theorem}

\begin{remark}
	It is worth noting that manifolds from items (i)--(iv) are actually manifolds from Theorem~\ref{thm:bert} but now the equators are not artificially removed from the manifold.
\end{remark}

Note that Theorem~\ref{thm:comp_bert_rim} is formulated for Riemannian manifolds only. Indeed, the pseudo-Riemannian case for manifolds with equators was not tackled before. In the present paper we prove that there are no completely Bertrand manifolds (with or without equators).

\begin{remark}
	We shall also note that~\cite{KudrFed3} also contains classification theorem for Bertrand manifolds with equators which uses the so-called {\em weakly Tannery surfaces} (defined and fully described in the same paper): 2-manifolds of revolution on which all {\em bounded} geodesics are closed. Study of the corresponding pseudo-Riemannian Bertrand manifolds lies beyond the scope of the present paper so we shall not formulate that theorem in fullness. We shall mention this topic again in the listing of problems for further research, see below.
\end{remark}

Now let us provide several theorems for Tannery manifolds.

Complete description of Riemannian Tannery manifolds can be found in~\cite{Besse}:

\begin{theorem}[\cite{Besse}]
	A spheroidal two-dimensional manifold $S\simeq (a,b)\times S^1$ of revolution is a Tannery manifold if and only if there exists such coordinates $(\psi, \varphi)$ where $\psi\colon (a,b) \to (0\pi)$ is a diffeomorphism, that the metric i those coordinates takes the form $$ds^2=R^2\left(\frac{1}{\beta^2}(1+h(\cos(\psi)))^2 d\psi^2 + \sin^2{\psi} d\varphi^2\right),$$ where $h\colon (-1,1)\to (-1,1)$ is some odd function, and $R\in \mathbb{R}$, $\beta \in \mathbb{Q}$ are positive constants.
\end{theorem}

A natural problem is to describe Tannery manifolds in pseudo-Riemannian case. Surprisingly few facts are known regarding this problem. The following interesting results are due to Mounoud and Mounoud-Suhr~\cite{Mon, MonSuhr}.

\begin{theorem}[\cite{Mon}]
	Consider a Lorentzian surface such that it is orientable, has a periodic spacelike Killing field, and all its geodesics are closed except possibly for the one orthogonal to the Killing field. Then it is projectively equivalent to an open subset of a Riemannian Tannery surface.
	
	Moreover, any Riemannian Tannery surface contains an open subset projectively equivalent to a Lorentzian surface with the described properties.
\end{theorem}

Delving further, let us say, following~\cite{MonSuhr}, that a pseudo-Riemannian metric is a {\em timelike (spacelike) SC-metric} if all timelike (spacelike) geodesics are simply closed with the same pseudo-Riemannian length (that is, an analogue of Zoll surfaces in Riemannian case).

\begin{theorem}[\cite{MonSuhr}] \label{thm:tannery}
	If $(M,g)$ is a Lorentzian 2-manifold all of whose timelike (spacelike) geodesics are closed, hen $(M,g)$ is finitely covered by a timelike (spacelike) SC-metric of $S^1\times \mathbb{R}$.
\end{theorem}

We shall use this theorem when discussing pseudo-Riemannian completely Bertrand manifolds.

\section{Certain results on pseudo-Riemannian Bertrand problem}

Theorem~\ref{thm:comp_bert_rim} together with the table from Theorem~\ref{thm:bert} shows an interesting fact: all Riemannian completely Bertrand manifolds with non-constant potentials are such, that the corresponding interval for $\theta$ in pseudo-Riemannian case is empty. This consideration makes the following question quite natural: are there any pseudo-Riemannian (to be precise, Lorentzian) Bertrand manifolds with non-constant potential? The following theorem gives answer to this (and even slightly broader) question.

\begin{theorem} \label{thm:no_completely}
	The are no completely Bertrand pseudo-Riemannian systems.
\end{theorem}

\begin{remark}
	Let us emphasise that in this context ``completely Bertrand system'' means that {\em all} trajectories are closed. If we require only spacelike (or timelike) trajectories to be closed, such systems do exist.
\end{remark}

\begin{proof}
	To prove this theorem we use a technique similar to the one used in \cite{KudrFed1, KudrFed2, KudrFed3} to tackle the Riemannian case.
	
	Suppose there exists a completely Bertrand pseudo-Riemannian system $(S, ds^2, V)$. The {\em Maupertuis principle} states (see, for example,~\cite{Chanda}) that for any sufficiently large energy level $E > \inf{V}$ all local solutions $x(t)$ of the equations of motion such that the energy value on them equals $E$ coincide (up to a reparametrization) with geodesics on the submanifold $S_E := V^{-1}(-\infty, E) \times S^1 \subseteq S$ with the metric $g_E=(E-V(x))ds^2$.
	
	Let us call this mapping of $(S,ds^2)$ to $(S_E, g_E)$ the {\em Maupertuis mapping}. Since the property ``to be closed'' of trajectories is preserved by Maupertuis mapping, if the original system was completely Bertrand, the resulting manifold $(S_E, g_E)$ would be a (Lorentzian) Tannery manifold in the sense that {\em all} geodesics on it are closed.
	
	But it follows from Theorem~\ref{thm:tannery} (see \cite[Theorem~6.4]{MonSuhr}) that every Lorentzian surface contains either a nonclosed timelike or a nonclosed spacelike geodesic.
	
	This contradiction completes the proof.
\end{proof}

Note that in the light of this theorem it is natural to study spacelike, timelike completely Bertrand manifolds. Their families are known to non-empty. To be precise, consider the 2-dimensional pseudosphere (or a family thereof parametrised by the radius $r>0$): $$S^2_1(r)=\{x\in\mathbb{R}^3 \;|\; \langle x,x \rangle_1 = -(x^1)^2+(x^2)^2=r^2 \}.$$ All spacelike geodesics on such manifold are closed, so a pseudosphere provides an example of a spacelike Tannery surface, and hence a spacelike completely Bertrand manifold.

Moreover, inverting all signs, this example gives pseudo-Riemannian manifolds on which all {\em timelike} geodesics are closed (and hence an example of timelike completely Bertrand system). \newline

Now let us move from the very restricting class of completely Bertrand manifolds to, on the contrary, very general class of locally Bertrand manifolds.

Let $S$ be a pseudo-Riemannian locally Bertrand manifold with the metric $ds^2=dr^2-f^2(r)d\varphi^2$ (possibly with equators). The definition of locally Bertrand system there exists at least one strongly stable circular orbit. This orbit may either be a) an equator or b) not an equator.

Let us consider the case when there exist strongly stable circular orbits of type b), that is, non-equatorial strongly stable orbits. If $S$ does not have any equators, by Theorem~\ref{thm:bert}, it is a belt of one of the manifolds listed above with the corresponding potential on it.

Therefore the interesting case is that of the manifold $S$ having equators. Let us remove those equators from $S$. The manifold is thus split into several pieces $(S_i, ds^2|_{S_i})$ with potentials $V|_{S_i}$. Note that the boundary of the closure of each piece $S_i$ is either an equator, or the boundary circle of the whole manifold $S$.

\begin{remark}
	The manifold $S$ is an open manifold diffeomorphic to $(a,b)\times S^1$. Saying ``the boundary circle'' of this manifold, we mean the circles $\{a\}\times S^1$ and $\{b\}\times S^1$. Note that the profile function $f(r)$ may go to zero or to infinity when $r$ goes to $a$ or to $b$. In that case the corresponding boundary circle degenerates into a point or a ``circle'' of infinite radius.
\end{remark}

Furthermore, for each $i$ either of the two following options takes place: either both connected components of $\partial \bar{S_i}$ are equators, or one component is an equator and the other one --- a boundary circle of $S$. Let us call the pieces of the first type ``interior'', and the ones of the second type --- ``exterior''

Consider a piece $S_i$ such that it contains a strongly stable circular orbit. Since the triple $(S,ds^2,V)$ was locally Bertrand, the triple $(S_i, ds^2|_{S_i})$ is locally Bertrand as well. It does not have any equators, so it is a belt of one of the manifolds listed in Theorem~\ref{thm:bert}. The explicit description of the geometry of these manifolds (see~\cite{Zagr}) shows that {\em i)} none of these manifolds have two equators as their boundary circles, {\em ii)} only the manifolds, corresponding to the parameters $c<0, \delta=0$ or $\delta<0, c<0, c^2+4\delta>0$ have an equator as a boundary circle.

That means that no interior piece may be a (locally) Bertrand system. Therefore, only an exterior piece may contain strongly stable circular orbits. Moreover, each exterior piece is a belt of a pseudo-Riemannian (locally) Bertrand manifold of the type 2 or 4 (as listed in Theorem~\ref{thm:bert}(3)). Hence we have the following theorem:

\begin{theorem}
	Let $(S,ds^2,V)$ be a locally Bertrand system. Then it is of the form ``cap -- black box -- cap'', where cap and the black box are glued along equators (the caps may be absent and the black box in that case may be cut away from the equators), each cap is a belt adjacent to an equator of a pseudo-Riemannian locally Bertrand manifold of type 2 or 4 from Theorem~\ref{thm:bert}(3) with the corresponding potential on it. The black box does not contain non-equatorial strongly stable circular orbits.
\end{theorem}

\section{Problems for further research}

In this section we formulate several problems regarding the theory of pseudo-Riemannian Bertrand systems for further research with brief commentary.

\begin{problem}
	Study spacelike and timelike pseudo-Riemannian completely Bertrand manifolds.
\end{problem}

As it was mentioned before, these classes of Betrand systems are non-empty because there exist spacelike and timelike Tannery manifolds. It seems reasonable to approach this problem, like in the Riemannian case, via Maupertuis principle. That leads to the following problem, which is interesting and important all by itself:

\begin{problem}
	Describe spacelike and timelike pseudo-Riemannian Tannery surfaces.
\end{problem}

The next interesting class of Bertrand systems is the systems with closing potential, that is, the systems all {\em bounded} trajectories of which are closed. Note that the arguments of Theorem~\ref{thm:no_completely} do not work in this case: it may be so that all trajectories which are mapped to timelike (or spacelike) geodesics via the Maupertuis mapping are not bounded and hence non-closed.

\begin{problem}
	Describe all pseudo-Riemannian Bertrand systems (systems with closing potential) and, in particular, weakly Tannery manifolds.
\end{problem}

As a final problem of this list, let us mention that in Riemannian case a complete diagram of inclusion of different classes of Bertrand systems was provided in~\cite{KudrFed3}. Little is know in that regard for pseudo-Riemannian systems. As it was mentioned, in the case of manifolds without equators it was proved (see~\cite{Zagr}) that the classes of Bertrand, locally, semilocally, strongly and weakly Bertrand systems coincide. It is important to obtain a similar diagram of classes in the general case.

\begin{problem}
	Describe the relationship between different classes of pseudo-Riemannian Bertrand systems in the form of a (strict) diagram of inclusions.	
\end{problem}

\begin{ack}
The work was supported by Russian Science Foundation, project 21-11-00355, and was done at Lomonosov Moscow State University.
\end{ack}

\end{document}